\documentclass{amsart}
\usepackage[latin1]{inputenc}
\usepackage[all]{xypic}

\advance\voffset 0cm
\advance\textheight 0cm
\advance\hoffset 0cm
\advance\textwidth 0cm

\newtheorem{Def}{Definition}[section]

\newtheorem{Lemma}[Def]{Lemma}
\newtheorem{Thm}[Def]{Theorem}
\newtheorem{Cor}[Def]{Corollary}

\def \Spec {{\rm Spec \, }}
\def \Hom {{\rm Hom \, }}

\def \jacobi #1 #2 {\left( \frac{\partial #1_i}{\partial #2_j} \right)_{i,j}}

\newcommand{\absvalue}[1]{\left| #1 \right|}

\begin{document}

 \title{Completeness and Compactness for Varieties over a local Field}
 \author{Oliver Lorscheid}
 \address{Department of Mathematics, Utrecht University , Budapestlaan 6, 3584 CD Utrecht, The Netherlands}
 \email{lorscheid@math.uu.nl}
 \thanks{The author likes to thank Florian Pop, Jakob Stix, Stefan Wewers, Gunther Cornelissen and his own parents for their support}
 \date{August 16, 2006}

 \begin{abstract}
  This paper shows that for a local field $K$, a subfield $k \subset K$ and a variety $X$ over $k$, $X$ is complete if and only if
  for every finite field extension $K' \mid K$, the set $X(K')$ is compact in its strong topology.
 \end{abstract}

 \maketitle

 \section{Introduction}

 The main statement Theorem (\ref{mainthm}) of this article is a fact well-known to experts, but which apparently misses a proof in the 
 literature. The complex case is proven in \mbox{\cite[I §10]{MUM99}}. For non-archimedean local fields, one can study rigid analysis 
 as introduced by Tate and further developed by Remmert, Gerritzen, Kiehl, and others, which perhaps might lead to an analytic proof 
 of the theorem. In the present paper however, we do not use methods from non-archimedean analysis, but give instead an algebraic proof.

 \medskip
 We define a local field to be a locally compact topological field. A topological field is a field with a non-discrete topology
 such that the field operations are continuous and points are closed.
 In \cite[I §2-4]{WEIL67} these fields were shown to be isomorphic and homeomorphic to the finite extensions of $\mathbb R$, $\mathbb Q_p$ or
 $\mathbb F_p((X))$ for some prime $p$.

 \begin{Thm}
  \label{mainthm}
  Let $K$ be a local field, $k \subset K$ a subfield and $X$ a variety \mbox{over $k$}. Then $X$ is complete if and only if 
  for every finite extension $K' \mid K$ of fields $X(K')$ is compact in its strong topology.
 \end{Thm}

 A pre-variety $X$ over a field $k$ is a geometrically reduced scheme of finite type over $k$, not necessarily irreducible.
 For a field extension \mbox{$K \mid k$}, the $K$-rational points of $X$ are defined to be the set
 $X(K) = \Hom_k (\Spec K, X) = \Hom_K (\Spec K, X \otimes_k K)$
 and can be viewed as a subset of the closed points of $X \otimes_k K$.

 For every topological field $k$, we endow the sets of $k$-rational points of pre-varieties over $k$ with a topology by requiring the following
 properties: 

 \begin{enumerate}
   \item $k \to \mathbb A^1_k(k)$ is a homeomorphism,
   \item $(X \times Y)(k) = X(k) \times Y(k)$ has the product topology of $X(k)$ and $Y(k)$, and
   \item for $Z \hookrightarrow X$ locally closed, $Z(k)$ has the subspace topology of $X(k)$.
 \end{enumerate}

 There is a unique way of defining a topological structure like this and we will call it the strong topology 
 on the pre-varieties over $k$. We remark that the strong topology is stronger than the Zariski topology. 
 In the following, properties of spaces or maps in the strong topology are marked with the adverb ``strongly'',
 such as ``strongly closed'' for ``closed in the strong topology''.

 If we have a morphism $\varphi: X \to Y$ of prevarieties over $k$, then $\varphi(k): X(k) \to Y(k)$ becomes strongly continuous. Similarly,
 if we have a pre-variety $X$ over $k$ and a homomorphism $k \to k'$ of topological fields which is a topological embedding, 
 then the induced map $X(k) \to X(k')$ also is a topological embedding. This means that properties like open, closed or dense in the
 strong topology are stable under such field extensions. In particular, extensions of fields with absolute value are of this form.
 By definition of a variety, the $k$-rational points of a $k$-variety form a Hausdorff space.
  
 \medskip
 To see that it is necessary to consider finite field extensions in \mbox{theorem (\ref{mainthm})} look at the example
 $X = \Spec \mathbb R [X_1,X_2]/(X_1^2 + X_2^2-1)$.
 Here $X(\mathbb R) \simeq \mathbb S^1$ is compact, but $X(\mathbb C) \simeq \mathbb C^\times$ is not.
 Another example is to take an irreducible polynomial $f(X_1)$ of degree $n\geq 2$ over a local field $K$.
 Then $X = \Spec K[X_1,X_2]/f(X_1)$ has no $K$-rational points, hence $X(K)$ is compact. But for the splitting field $K'=K(f)$,
 $X(K')$ is the union of $n$ affine spaces which certainly is not compact.

 \bigskip

 \bigskip

\section{From Completeness to Compactness}

 \bigskip

 In this section, $K$ denotes a local field and $\absvalue \ $ its absolute value.

 The key result for the proof that a complete variety $X$ over $K$ is compact in its strong topology is Chow's Lemma.
 A proof can be found in \mbox{\cite[5.6.1]{EGAII}}, also compare with \mbox{\cite[I §10]{MUM99}}.

 \begin{Thm}
  \textrm{\bf (Chow´s Lemma)}

  Let $X$ be a reduced scheme that is proper and of finite type over a noetherian scheme $S$.
  Then there is a reduced projective scheme $Y$ of finite type over $S$,
  a surjective morphism $\varphi: Y \to X$ over $S$, and an open and dense subset $U \subset X$
  such that the restriction $\varphi: \varphi^{-1}(U) \to U$ is an isomorphism.
 \end{Thm}

 Using the projection of the $(n+1)$-sphere to $\mathbb P^n(K)$, one sees that the $K$-rational points of projective space 
 form a strongly compact set. For this, it is crucial that $K$ is locally compact.
 Therefore, strongly closed subsets of projective space are also strongly compact, and we obtain:

\begin{Lemma}
  For a projective variety $X$ over $K$, the set $X(K)$ is strongly compact.
\end{Lemma}

 Now we can conclude via Chow's Lemma:

\begin{Thm}
  More generally, for a complete variety $X$ over $K$, the set $X(K)$ is strongly compact.
\end{Thm}

\begin{proof}
  The proof is by induction with respect to $\dim X$.

  For $\dim X = 0$, $X(K)$ is a finite set and therefore compact.

  For $\dim X > 0$, by Chow's Lemma there is a projective variety $Y$ and a surjective morphism $\varphi: Y \to X$ as well as
  an open and dense subvariety $U \subset X$ such that $\varphi: V \to U$ is an isomorphism for $V = \varphi^{-1}(U)$.

  We obtain the following diagram of sets and maps:

  $$
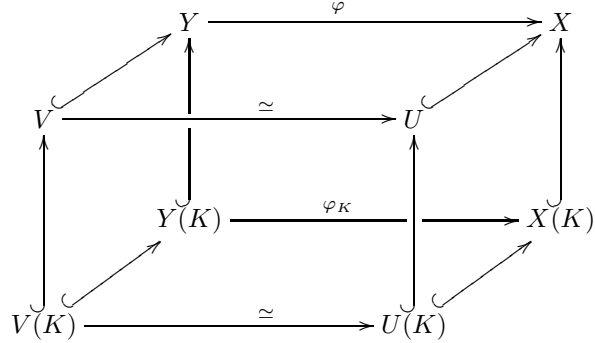
\xymatrix{
                                         & Y \ar[rrr]^(.4){\varphi} &             &                    & X                     \\
  V \ar@{^{(}->}[ur] \ar[rrr]^(.6)\simeq &                          &             & U \ar@{^{(}->}[ur] &                       \\
                                         & Y(K) \ar[rrr]^(.4){\varphi_K}|(.605)\hole \ar@{^{(}->}[uu]|\hole
                                                                    &             &                    & X(K) \ar@{^{(}->}[uu] \\
  V(K) \ar@{^{(}->}[ur] \ar[rrr]^(.6)\simeq \ar@{^{(}->}[uu] &      &             & U(K) \ar@{^{(}->}[ur] \ar@{^{(}->}[uu]      }$$

  \bigskip

  As an isomorphism between $V$ and $U$, $\varphi$ is a bijection between their $K$-rational points
  and \mbox{$U(K) = U \cap X(K)$} is part of the image of $\varphi_K: Y(K) \to X(K)$.

  Now the closed subvariety $Z = X - U$ has the same $K$-rational points as $X$ outside of $U$. 
  Since $U \subset X$ is an open and dense subvariety, $\dim Z < \dim X$, and since $Z \to X$ is a closed subscheme, $Z$ is complete. 
  Hence we can assume by induction that $Z(K)$ is strongly compact.

  Finally, as the union of the two strongly compact sets $Z(K)$ and $\varphi_K(Y(K))$, $X(K)$ is itself strongly compact.
 \end{proof}

 Now $X(K) = Hom_K(K,X \otimes_k K)$ if $X$ is defined over a subfield $k \subset K$.
 Since completeness is stable under base extension, and keeping in mind that we can replace $K$ by any finite extension $K'$,
 the first half of theorem (\ref{mainthm}) follows:

 \begin{Cor}
  If $X$ is a complete variety over $k \subset K$, then $X(K)$ is strongly compact.
 \end{Cor}

 \medskip
 \bigskip

\section{The Implicit Function Theorem}

 \bigskip

 To prove the inverse implication, we need an analogon of the implicit function theorem for spaces of rational points over a local field
 if the map in question comes from a morphism of varieties. The theorem will provide only a local homeomorphism,
 which follows from Hensel's Lemma for the non-archimedean case.

 \begin{Thm}
  \textrm{\bf (Hensel´s Lemma)}

  Let $K$ be a complete field relating to a non-archimedean absolute value $\absvalue \ $, $R$ the ring of integers and
  $f_1, \ldots, f_n \in R[X_1,, \ldots X_n]$, $a_1, \ldots, a_n \in R$ with

  $$ \absvalue{f_l(a_1, \ldots, a_n)} < \absvalue{\det \jacobi f X (a_1, \ldots, a_n)}^2 $$

  \medskip

  for $l=1,\ldots,n$. Then there are unique $b_1, \ldots, b_n \in R$ with

  $$ f_l(b_1, \ldots, b_n) = 0 \hspace{.5cm} \textrm{and}\hspace{.5cm}\absvalue{b_l - a_l} < \absvalue{\det \jacobi f X (a_1, \ldots, a_n)} \ \ $$

  \medskip

  for $l=1,\ldots,n$.
 \end{Thm}

 Except from uniqueness one may recognize this to be a special case of the ``General Hensel Lemma'' in \cite[(5.21)]{GRE69}
 and uniqueness is easy to prove.

 \begin{Thm}
  \label{Thm: impl func}
  \textrm{\bf (Implicit function theorem)}

  Let $K$ be a complete field relating to an absolute value $\absvalue \ $, $X$ a $r$-dimensional variety over $K$
  and $\varphi: X \to \mathbb{A}_K^r$ a morphism, that is étale at a  $K$-rational point $p \in X$ and maps $p$ to a $K$-rational point.
  Then $\varphi_K: X(K) \to \mathbb{A}_K^r(K)$ is a homeomorphism strongly near $p$.
 \end{Thm}

 \begin{proof}
  If $K$ is archimedean, by Ostrowski's theorem (\cite[II 4.2]{NEU92}) $K = \mathbb R$ or $\mathbb C$,
  and the theorem follows from real analysis.

  Let $K$ be non-archimedean. For $\varphi$ to be étale at $p$ means that near $p$ it is defined by a ring homomorphism of the form

  $$ K[Y_1,\ldots,Y_r] \ \ \longrightarrow \ \ K[Y_1,\ldots,Y_r][X_1,\ldots,X_n] / (f_1,\ldots,f_n) $$

  $$ \textrm{with} \ \ \ \det \jacobi f X (p) \neq 0 \ \ \ \textrm{in} \ \ k(p) = K.$$

  \medskip

  Because both $p$ and $\varphi(p)$ are $K$-rational, we can assume by translation
  that $p$ has local coordinates $(0,\ldots0,a_1,\ldots,a_n)$ in $\mathbb{A}_K^{r+n}(K) \simeq \mathbb{A}_{K[Y_1,\ldots,Y_r]}^n(K)$
  with \mbox{$a_l \in R \subset K$}, the ring of integers. By multiplication with a sufficiently small number
  in $K$, we can assume that $f_1,\ldots,f_n \in R[Y_1,\ldots,Y_r][X_1,\ldots,X_n]$ without changing coordinates. Therefore

  $$ \delta \ := \ \absvalue{\det \jacobi f X (p)} \ = \ \absvalue{\det \jacobi f X (0,\ldots,0,a_1,\ldots,a_n)} \ \leq \ 1 . $$

  \medskip

  When $X_l = a_l$, there is no constant term of the $f_l$ in $Y_1,\ldots Y_r$.
  Hence for $y_1,\ldots,y_r \in K$ with $\absvalue{y_l} < \delta^2$, it follows that

  $$ \absvalue{f_l(y_1,\ldots,y_r,a_1,\ldots,a_n)} \ < \ \delta^2. $$

  \medskip

  By Hensel's Lemma there are unique $b_1,\ldots,b_n$ with
  $f_l(y_1,\ldots,y_r,b_1,\ldots,b_n) = 0$ and $\absvalue{b_l - a_l} < \delta$. This defines an inverse map
  $\varphi_K^{-1}: B_{\delta^2} (0) \to B_{\delta}(p)$, where $B_\epsilon (q) = \{ r \mid \absvalue{r-q} < \epsilon \}$ is an open neighbourhood
  in the appropriate space.

  To prove the continuity of $\varphi_K^{-1}$, consider that for $\absvalue{y_l} < \delta$

  $$ \absvalue{\det \jacobi f X (y_1,\ldots,y_r,a_1,\ldots,a_n)} \ = \delta $$

  \medskip

  is constant. Hence it is sufficient to show continuity at $0$ and this is proven
  once there are points arbitrarily close to $p$ in the image of the inverse map.

  Let $\epsilon > 0$, then multiply $f_1$ with a sufficient small number in $K$, such that $\delta \leq \epsilon$. Choose $y_1\ldots,y_r \in K$
  with $\absvalue{y_l} < \delta^2$ and find a point $(y_1,\ldots,y_r,b_1,\ldots,b_n)$ with $f_l(y_1,\ldots,y_r,b_1,\ldots,b_n) = 0$ and
  $\absvalue{b_l - a_l} < \delta \leq \epsilon$ by Hensel's Lemma.
 \end{proof}

 \bigskip

 \bigskip

\section{From Compactness to Completeness}

 \bigskip

 This part of the proof holds for any field $K$ complete with respect to an absolute value. Let $k$ be a subfield of $K$.
 We will need the following theorem.  A proof can be found in \cite{NAG62}. 

 \begin{Thm}
  \label{nagata}
  \textrm{\bf (Nagata)}

  Let $k$ be an arbitrary field.
  Every variety over $k$ can be embedded into a complete variety over $k$ as an open and dense subset.
 \end{Thm}

 \medskip
 Now it is possible to prove the other direction of theorem (\ref{mainthm}).

 \begin{Thm}
  Let $X$ be a variety over $k$ such that for every finite extension $K' \mid K$, $X(K')$ is strongly compact.
  Then $X$ is complete.
 \end{Thm}

 \begin{proof}
  If $\dim X = 0$, there is nothing to prove. Let $\dim X > 0$.

  By Nagata's theorem, $X$ can be embedded into a complete variety $Y$ over $k$ as an open and dense subset.
  If $X$ was not complete there would be a closed point $p \in Y - X$.
  But we will show that the existence of $p$ leads to a contradiction.

  Since $X \neq Y$ implies $X \otimes k' \neq Y \otimes k'$ for every extension $k' \mid k$, we can change the base field
  to prove the non-existence of the point $p$. Further, we can substitute $K$ by a finite extension $K'$
  without changing the validity of the conditions of the theorem.
  First, we assume that $X$ and $Y$ are varieties over $K$ and that $p$ is $K$-rational.

  In an affine neighbourhood of $p$, we can find a curve through $p$ that intersects $X$ non-trivially. Let $C$ be its closure in $Y$,
  then $C \cap X$ is closed in $X$ and therefore $(C \cap X)(K')$ strongly compact for every $K' \mid K$ finite.
  Replace $Y$ by $C$ and $X$ by $C \cap X$, then $Y - X$ has only finitely many points.

  Since $X(K)$ is strongly compact, it is strongly closed in $Y(K)$. It follows that $p$ is a strongly discrete point as the complement
  of finitely many strongly closed points in the strongly open set $(Y-X)(K)$.

  The normalisation $Z \to Y$ is a finite morphism, and therefore there are finitely many discrete points in the inverse image of $p$.
  Since the normalisation of a curve over an algebraically closed field is smooth, we can enlarge $K$ such that the normalisation of $Y$ 
  is smooth. Hence we can assume that $Y$ is a normal and smooth curve.

  A smooth curve is locally étale over $\mathbb{A}_K^1$, compare with \cite[III 10.4]{HAR77} and \cite[III §6 Thm.1]{MUM99}.
  This means that there is an open neighbourhood $U \subset Y$ of $p$ and an étale morphism $\varphi: U \to \mathbb A_K^1$.
  By another finite field extension, we can assume that not only $p$, but also $\varphi(p)$ is $K$-rational.
  Then the implicit function theorem (\ref{Thm: impl func}) shows that there are homeomorphic strong neighbourhoods of $p$
  and $\varphi(p)$. But $\varphi(p)$ cannot be strongly discrete since $\mathbb A_K^1(K) \cong K$ is non-discrete. Contradiction!
 \end{proof}


\bigskip

 \bigskip

 \begin{thebibliography}{AutorJahr}
  \bibitem[1]{EGAII}   A. Grothendieck, J. Dieudonné, {\em Eléments de Géométrie \mbox{Algébrique II}}, Publ. Math. \mbox{IHES 8}, 1961
  \bibitem[2]{GRE69}   Marvin J. Greenberg, {\em Lectures on Forms in Many Variables}, W. A. Benjamin Inc., New York-Amsterdam 1969
  \bibitem[3]{HAR77}   Robin Hartshorne, {\em Algebraic Geometry}, Springer-Verlag, New York 1977
  \bibitem[4]{NAG62}   Nagata, {\em Imbedding of an Abstract Variety in a Complete Variety},
                            Journal for Mathematics of the Kyoto University 2, 1962, pages 1-10
  \bibitem[5]{NEU92}   Jürgen Neukirch, {\em Algebraische Zahlentheorie}, Springer-Verlag, Berlin-Heidelberg 1992
  \bibitem[6]{MUM99}   David Mumford, {\em The Red Book of Varieties and Schemes}, second edition, Springer-Verlag, Berlin-Heidelberg 1999
  \bibitem[7]{WEIL67} André Weil, {\em Basic Number Theory}, Springer-Verlag, Berlin-Heidelberg-New York 1967
 \end{thebibliography}
\end{document}